\documentclass[a4,12pt,twoside]{article}
\usepackage{amsfonts}
\usepackage{amssymb}
\usepackage{amsthm}
\usepackage{amsmath}
\usepackage[dvips]{graphicx,color}
\usepackage{latexsym,bm}

\numberwithin{equation}{section}

\def\R2{\mathbb{R}^2}

\begin{document}           

\begin{flushleft}
\renewcommand{\thefootnote}{\fnsymbol{footnote}}
{\bf\LARGE  Continuity of packing measure function of self-similar
iterated function systems
}  
\footnotetext{{\bf MR(2000) Subject Classification:}
\hspace{0.2cm}28A78, 28A80 }
 \footnotetext{Supported by NSFC 10901081}
 \vskip0.5cm {\large
Hua Qiu\hspace{0.1cm} }

 {\raggedright \small  Department of Mathematics,
Nanjing University, Nanjing, 210093, China}
\\ { Email: huatony@gmail.com}

\end{flushleft}

\begin{abstract}
In this paper, we focus on the packing measure of self-similar sets.
Let $K$ be a self-similar set whose Hausdorff dimension and packing
dimension equal $s$, we state that if $K$ satisfies the strong open
set condition with an open set $\mathcal{O}$, then
 $$\mathcal{P}^s(K\cap B(x,r))\geq (2r)^s$$ for each open ball $B(x,r)\subset \mathcal{O}$ centered in $K$, where $\mathcal{P}^s$ denotes the $s$-dimensional packing measure. We use this inequality to obtain some precise density theorems for packing measure of self-similar sets, which can be applied to compute the exact value of the $s$-dimensional packing measure $\mathcal{P}^s(K)$ of $K$. Moreover, by  using  the above results, we show the continuity of the packing measure function of the attractors on the space of self-similar iterated function systems satisfying the strong separation condition. This result gives a complete answer to a question posed by L. Olsen in $\cite{Ols}$.
\end{abstract}

\section{Introduction}
In this paper we will analysis the behaviour of the packing measure
of self-similar sets with open set condition or strong separation
condition. Recall the definition of packing measure, introduced by
Tricot \cite{Tri}, Taylor and Tricot \cite{TT}, which requires two
limiting procedures. For $E\subset\mathbb{R}^d$ and $\delta>0$, a
$\delta$-packing of $E$ is a countable family of disjoint open balls
of radii at most $\delta$ and with centers in $E$. For $s\geq 0$,
the $s$-dimensional \textit{packing premeasure} of $E$ is defined as
$$P^s(E)=\inf_{\delta>0}\{P^s_\delta(E)\},$$
where $P^s_\delta(E)=\sup\{\sum_{i=1}^{\infty} \mbox{diam}(B_i)^s\}$
with the supremum taken over all $\delta$-packing of $E$. Here
$\mbox{diam}(B_i) $ denotes the diameter of $B_i$. The
$s$-dimensional \textit{packing measure} of $E$ is defined as
$$\mathcal{P}^s(E)=\inf\{\sum_{i=1}^\infty P^s(E_i)| E\subset\bigcup_{i=1}^\infty E_i\}.$$
The \textit{packing dimension} of $E$ is defined as
$$\dim_{P}(E)=\inf\{s\geq 0| \mathcal{P}^s(E)=0\}=\sup\{s\geq 0| \mathcal{P}^s(E)=\infty\}.$$ The packing measure and packing dimension play an important role in the study of fractal geometry in a manner dual to the Hausdorff measure and Hausdorff dimension (See \cite{Fal} and \cite{Mat1} for further properties of the above measures and dimensions).

Let $N\geq 2$ be an integer. Let $\textbf{f}=\{f_1,f_2,\cdots,f_N\}$
be an \textit{iterated function system} (\emph{IFS}) on
$\mathbb{R}^d$
 of contractive similitudes. The corresponding \textit{self-similar set} for $\textbf{f}$ is the unique non-empty compact set $K\subset \mathbb{R}^d$ which is invariant under the action of the elements of $\textbf{f}$:
$$K=\bigcup_{i=1}^N f_i(K).$$
It is well-known that if  $\textbf{f}$ satisfies the \textit{open
set condition}(\emph{OSC}), i.e., there exists a nonempty bounded
open set $\mathcal{O}\subset \mathbb{R}^d$ such that
$f_i(\mathcal{O})\subset \mathcal{O}$ for all $1\leq i\leq N$ and
$f_i(\mathcal{O})\cap f_j(\mathcal{O})=\emptyset$ for all  $i\neq
j$, then the Hausdorff dimension $\dim_H(K)$ and the packing
dimension $\dim_P(K)$ of $K$ coincide, and the common value
$s=\dim_H(K)=\dim_P(K)$ is given by the following formula
\begin{equation}\label{0}
\sum_{i=1}^N r_i^s=1,
\end{equation}
where $r_i$ denotes the contraction ratio of $f_i$ for $1\leq i\leq
N.$ Moreover, the Hausdorff measure and packing measure of $K$ are
finite and positive. This was proved by Moran ${\cite{Mor}}$ in 1946
and rediscovered by Hutchinson ${\cite{Hut}}$ in the 1980s. Since
the intersection of $\mathcal{O}$ and $K$ may be empty, the OSC is
in general too weak to imply results. One can strengthen the
definition as follows: \textit{The strong open set condition (SOSC)}
holds if and only if furthermore $\mathcal{O}\cap K\neq \emptyset.$
Schief proved that SOSC is equivalent to OSC in the Euclidean case,
see $\cite{Sch}$. There is another separation condition called the
\textit{strong separation condition (SSC)} which is satisfied if
$f_i(K)\cap f_j(K)=\emptyset$ for all $i,j$ with $i\neq j.$
Obviously, SSC implies SOSC and the implication may not be inverted.
In this paper, we will frequently assume these two conditions.

We shall need some standard notations from symbolic dynamics. For
each positive integer $k$, let
$$W_k=\{1,2,\cdots,N\}^k=\{\textbf{i}=i_1i_2\cdots,i_k: i_j\in
\{1,2,\cdots, N\}\},$$ denote the space of words of length $k$ with
symbols $\{1,2,\cdots,N\}$. Also, for $k=0$, we define
$W_0=\{\emptyset\}$ and call $\emptyset$ the empty word. Moreover,
set $W=\bigcup_{k\geq 0}W_k$ and denote the length of $\textbf{i}\in
W$ by $|\textbf{i}|$. Assume now that $\textbf{f}=\{f_1,f_2\cdots,
f_N\}$ is an \textit{IFS} with invariant set $K$. Let
$f_{\emptyset}=id, r_{\emptyset}=1, K_{\emptyset}=K.$ For each
non-empty finite word $\textbf{i}=i_1i_2\cdots i_k$ let
$f_\textbf{i}=f_{i_1}\circ f_{i_2}\circ\cdots\circ f_{i_k}$,
$r_\textbf{i}=r_{i_1}r_{i_2}\cdots r_{i_k}$ and
$K_\textbf{i}=f_\textbf{i}(K)$. Then $K=\bigcup_{\textbf{i}\in
W_k}K_\textbf{i}$ for each $k$.

For a Borel measure $\mu$ on $\mathbb{R}^d$ and a Borel set $E$, we
let $\mu|_E$ denote the restriction of $\mu$ to $E$. Let $\lambda$
denote the self-similar measure satisfying
$$
\lambda=\sum_i r_i^s\lambda\circ f_i^{-1}.
$$
It is well known that under the assumption of OSC,
$\lambda=\frac{\mathcal{H}^s|_K}{\mathcal{H}^s(K)}=\frac{\mathcal{P}^s|_K}{\mathcal{P}^s(K)}$,
and $\lambda(K_\textbf{i})=r_\textbf{i}^s$ for each $\textbf{i}\in
W.$ Hence the measure $\lambda$ is the normalised $s$-dimensional
Hausdorff measure restricted to $K$.

We always assume that $K$ is in general position, i.e., not
contained in a hyperplane. In \cite{Mat} it is proved that, under
this assumption and the OSC, the intersection of $K$ with any
$l$-dimensional $C^1$ submanifold of $\mathbb{R}^d$ with $0<l<d$ is
an $\mathcal{H}^s$- null set, and therefore a $\mathcal{P}^s$- null
set.

Since the definitions of Hausdorff and packing measures are
sometimes awkward to work with, there are only very few non-trivial
examples of sets in $E\subset\mathbb{R}^d$ for which the exact
Hausdorff measure $\mathcal{H}^{\dim_H(E)}(E)$ or packing measure
$\mathcal{P}^{\dim_P(E)}(E)$ of $E$ is known. For example, one can
see papers $\cite{AS}, \cite{Fen2}$. $\cite{ZF}$ is a recent review
of relevant open questions in this field. In particular, there is no
formula similar to $(\ref{0})$ for the Hausdorff measure or packing
measure of a self-similar set. In view of this, it is natural to ask
if the Hausdorff measure and packing measure vary continuously with
the IFS.

To make the above question precise we introduce the following
notations.

Let $N$ be a positive integer with $N\geq 2$ and let $X\subset
\mathbb{R}^d$ be a compact set. Let
$\textbf{f}=(f_1,f_2,\cdots,f_N)$ be a IFS on $X$ satisfying the
OSC. In order to emphasize the relation between the corresponding
fractal characteristics with $\textbf{f}$,  we write $K(\textbf{f})$
for the self-similar set associated with  $\textbf{f}$ and we write
$s(\textbf{f})$ for the common value of the Hausdorff dimension and
the packing dimension of $K(\textbf{f})$. Also, let
$\lambda(\textbf{f})$ denote the normalised
$s(\textbf{f})$-dimensional Hausdorff measure restricted to
$K(\textbf{f})$, and write the contraction ratio of $f_i$ as
$r_{i}(\textbf{f})$ for each $1\leq i\leq N$ and
$r_*(\textbf{f})=\min_{1\leq i\leq N}r_i(\textbf{f})$.

Write
$$M_{OSC}=\{\textbf{f}| \textbf{ f}=(f_1,f_2,\cdots,f_N)\mbox{ is a self-similar IFS on } X \mbox{  satisfying the OSC}\},$$
$$M_{SSC}=\{\textbf{f}|\textbf{ f}=(f_1,f_2,\cdots,f_N)\mbox{ is a self-similar IFS on } X \mbox{  satisfying the SSC}\}.$$
It is obvious that $M_{SSC}\subset M_{OSC}$. We equip $M_{OSC}$ and
$M_{SSC}$ with the metric induced by
$$D(\textbf{f},\textbf{g})=\max_{1\leq i\leq N}\{\parallel f_i-g_i\parallel_\infty\},$$ for $\textbf{f},\textbf{g}\in M_{OSC}$.
 It is not difficult to see
that $M_{SSC}$ is an open subset of $M_{OSC}$. Below
$d(\cdot,\cdot)$ denotes the Euclidean metric between two points or
two sets. For $\Delta>0$ we write
$$M_\Delta=\{\textbf{f}|\textbf{ f}\in M_{SSC},\mbox{ and } d(K_i(\textbf{f}),K_j(\textbf{f}))>\Delta \mbox{ for all } i\neq j\}.$$ Also, one can easily prove that
$M_{SSC}=\bigcup_{\Delta>0}M_{\Delta}$ and each $M_{\Delta}$ is an
open subset in $M_{SSC}$, see \cite{Ols}. The metric spaces
$M_{OSC}$ and $M_{SSC}$ provide a natural setting for investigating
to what extent fractal characteristics of $K(\textbf{f})$ vary
continuously with $\textbf{f}$. For example, let $\mathcal{K}$
denote the family of non-empty compact subsets of $X$ equipped with
the Hausdorff metric, then the map
$$\textbf{f}\rightarrow K(\textbf{f})$$ from $M_{OSC}$ into $\mathcal{K}$ is continuous, see \cite{CV}.
It also follows immediately from $(\ref{0})$ that the dimension map
$$\textbf{f}\rightarrow s(\textbf{f})$$ from $M_{OSC}$ into $\mathbb{R}$ is continuous. Following this line of investigation it is natural to ask if the measure maps
\begin{equation}\label{01}
\textbf{f}\rightarrow \mathcal{H}^{s(\textbf{f})}(K(\textbf{f})),
\end{equation}
\begin{equation}\label{02}
\textbf{f}\rightarrow \mathcal{P}^{s(\textbf{f})}(K(\textbf{f}))
\end{equation} from $M_{OSC}$ into $\mathbb{R}$ are continuous.

For the Hausdorff measure map, Ayer and Strichartz $\cite{AS}$
showed that in the special case of linear Cantor sets they found a
point $\textbf{f}\in M_{OSC}\setminus M_{SSC}$, at which the map
$(\ref{01})$ fails to be continuous. In $\cite{Ols}$ Olsen altered
the space $M_{OSC}$ to $M_{SSC}$, then he proved the continuity,
i.e., the map
$$\textbf{f}\rightarrow \mathcal{H}^{s(\textbf{f})}(K(\textbf{f}))$$ from $M_{SSC}$ into $\mathbb{R}$ is continuous.

Recall that in order to prove this continuity theorem, Olsen used a
so-called explicit formula for the Hausdorff measure of self-similar
sets which was established in $\cite{Ols1}$. Indeed, he showed that
the Hausdorff measure coincides with the infimum of the reciprocal
densities. Let  $\textbf{f}\in M_{OSC}$ and $K(\textbf{f})$ be the
corresponding self-similar set. Let $s(\textbf{f})$ be the Hausdorff
dimension and $\lambda(\textbf{f})$ be the normalised
$s(\textbf{f})$-dimensional Hausdorff measure restricted to
$K(\textbf{f})$. Then
$$\mathcal{H}^{s(\textbf{f})}(K(\textbf{f}))=\inf\{\frac{\mbox{diam}(U)^{s(\textbf{f})}}{\lambda(\textbf{f})(U)}| U \mbox{ is open and convex }, U\cap K(\textbf{f})\neq \emptyset\}.$$
Moreover, if furthermore $\textbf{f}$ satisfies the SSC, i.e., there
exists $\Delta>0$ such that $\textbf{f}\in M_\Delta$, then
\begin{equation}\label{03}
\mathcal{H}^{s(\textbf{f})}(K(\textbf{f}))=\inf\{\frac{\mbox{diam}(U)^{s(\textbf{f})}}{\lambda(\textbf{f})(U)}|
U \mbox{ is open and convex }, U\cap K(\textbf{f})\neq \emptyset,
\mbox{ diam}(U)\geq\Delta\}.
\end{equation}
As pointed in $\cite{Ols1}$, the above formulae are implicit in
earlier work by Marion and Ayer \& Strichartz, see
\cite{Mar},\cite{Mar1} and \cite{AS}. They used these formulae to
compute the exact value of the $s$-dimensional Hausdorff measure
$\mathcal{H}^s(\mathcal{C})$ of certain linear Cantor subsets
$\mathcal{C}$ of $\mathbb{R}$, where $s$ denotes the Hausdorff
dimension of $\mathcal{C}$.

By the fact that $M_{SSC}=\bigcup_{\Delta>0}M_{\Delta}$ and each
$M_{\Delta}$ is an open subset in $M_{SSC}$, the proof of Hausdorff
measure continuity theorem can be simplified to prove that  map
$(\ref{01})$ from $M_\Delta$ into $\mathbb{R}$  is continuous for
each $\Delta>0$. This is based on formula $(\ref{03})$ for the
Hausdorff measure $\mathcal{H}^{s(\textbf{f})}(K(\textbf{f}))$ of
$K(\textbf{f})$.  It should be pointed that  the fact that the
infimum in $(\ref{03})$ is taken over sets $U$ with $\mbox{
diam}(U)\geq\Delta$ plays a crucial role for the estimates involved
in the proof. Actually, A detailed comparison of the normalised
$s(\textbf{f})$-dimensional Hausdorff measure of elements in those
families where the infimum are taken over between one IFS and its
nearby IFSs is needed. The condition $\mbox{ diam}(U)\geq\Delta$
ensures the comparison successfully.

However, for the packing measure map, it is not clear if the map
$(\ref{02})$ is continuous. Let $\mathcal{A}$ denote the family of
analytic subsets of $M_{SSC}$ and $\sigma(\mathcal{A})$ denote the
$\sigma$-algebra generated by $\mathcal{A}$. Olsen  proved that the
map $(\ref{02})$ from $M_{SSC}$ into $\mathbb{R}$ is
$\sigma(\mathcal{A})$-measurable. However, this result is very weak
compared to the continuity, see  $\cite{Ols}$. The proof also needs
an explicit formula for the packing measure. Actually, a similar
formula for packing measure in the SSC case was also proved in
$\cite{Ols1}$, i.e., for each $\textbf{f}\in M_\Delta$ with
$\Delta>0$, the following formula holds.
\begin{equation}\label{04}
\mathcal{P}^{s(\textbf{f})}(K(\textbf{f}))=\sup\{\frac{(2r)^{s(\textbf{f})}}{\lambda(\textbf{f})(B(x,r))},
x\in K(\textbf{f}), 0<r\leq\frac{1}{2}\Delta\}.
\end{equation}
However, contrary to $(\ref{03})$, the explicit formula $(\ref{04})$
for the packing measure $\mathcal{P}^{s(\textbf{f})}$ of
$K(\textbf{f})$ does not involve taking infimum over sets whose
diameters are bounded away from zero. Indeed, in $(\ref{04})$ the
infimum is taken over all balls with radii less than
$\frac{1}{2}\Delta$. In particular, the infimum in $(\ref{04})$ is
taken over balls with arbitrarily small radii. For this reason it is
not possible to adapt the arguments in the proof of Hausdorff
measure continuity theorem (proving continuity of the
$\textbf{f}\rightarrow \mathcal{H}^{s(\textbf{f})}(K(\textbf{f}))$
using $(\ref{03})$) to prove continuity of the map
$\textbf{f}\rightarrow \mathcal{P}^{s(\textbf{f})}(K(\textbf{f}))$
using $(\ref{04})$.

On the other hand, results by Mattila \& Mauldin $\cite{Mau}$ show
that various somewhat related maps are discontinuous(in fact, not
even Borel measurable), and it is therefore entirely plausible that
the map $\textbf{f}\rightarrow
\mathcal{P}^{s(\textbf{f})}(K(\textbf{f}))$ is discontinuous.

Based on the above reasons,  Olsen posed the following open question
in $\cite{Ols}$.

\textbf{Question.} \textit{Is the packing measure function in
$(\ref{02})$ from $M_{SSC}$ into $\mathbb{R}$ continuous? If it is
not continuous, is it of Baire class $n$ for some positive integer
$n$? If it is not of Baire class $n$ for some positive integer $n$,
is it Borel measurable?}

Somewhat surprisingly, in this paper we will show that this map is
continuous, which gives a complete answer to this question. This
leads to our main result.

\textbf{Theorem 1.1.} \textit{The map $\textbf{f}\rightarrow
\mathcal{P}^{s(\textbf{f})}(K(\textbf{f}))$ from $M_{SSC}$ into
$\mathbb{R}$ is continuous.}

To prove Theorem 1.1, we first establish a new explicit formula for
the packing measure $\mathcal{P}^{s(\textbf{f})}(K(\textbf{f}))$
where the infimum is taken over all balls with radii bounded away
from zero. This formula is analogous to $(\ref{03})$(where the
infimum also is taken over sets with diameters bounded away from
zero) and replaces Olsen's formula $(\ref{04})$. The formula is
stated below. If $\textbf{f}\in M_\Delta$, $\Delta>0$, then
\begin{equation}\label{05}
\mathcal{P}^{s(\textbf{f})}(K(\textbf{f}))=\sup\{\frac{(2r)^{s(\textbf{f})}}{\lambda(\textbf{f})(B(x,r))},
x\in K(\textbf{f}), \frac{1}{2}r_*(\textbf{f})\Delta\leq
r\leq\frac{1}{2}\Delta\},
\end{equation}
where $r_*(\textbf{f})=\min_{1\leq i\leq N} r_i(\textbf{f}).$ Next,
using $(\ref{05})$, we then adapt the techniques for proving the
Hausdorff measure continuity theorem to establish Theorem 1.1
showing that the map $\textbf{f}\rightarrow
\mathcal{P}^{s(\textbf{f})}(K(\textbf{f}))$ is continuous. Here we
should point that the key formula $(\ref{05})$ is a direct corollary
of the following theorem which is another main result in this paper.

\textbf{Theorem 1.2.} \textit{Let $K\subset \mathbb{R}^d$ be the
self-similar set associated with a IFS $\textbf{f}=\{f_1, f_2,
\cdots, f_N\}$ satisfying the SOSC for an open set $\mathcal{O}$,
and with packing and Hausdorff dimension $s$. Then
\begin{equation}\label{1}
\mathcal{P}^s(K\cap B(x,r))\geq (2r)^s
\end{equation}
 for any ball $B(x,r)\subset\mathcal{O}$ centered in $K.$}

A similar result in the SSC case was proved by Olsen in
$\cite{Ols1}$, and from which the  explicit formula $(\ref{04})$ was
obtained.

It is not known whether the packing measure continuity theorem for
the map $\textbf{f}\rightarrow
\mathcal{P}^{s(\textbf{f})}(K(\textbf{f}))$ still holds from
$M_{OSC}$ into $\mathbb{R}$.  In the special setting of linear
Cantor sets of real line $\mathbb{R}$ Feng $\cite{Fen2}$ discussed
the exact value of packing measure
$\mathcal{P}^{\dim_P(\mathcal{C})}$ of self-similar Cantor sets
$\mathcal{C}$ satisfying the OSC(where the open set is an interval).
His result implies that $\mathcal{P}^{\dim_P(\mathcal{C})}$ depends
continuously on the IFSs in $M_{OSC}$. In view of this, we guess
that Theorem 1.1 could be generalized to the following setting.

\textbf{Conjecture.} \textit{The packing measure function in
$(\ref{02})$ from $M_{OSC}$ into $\mathbb{R}$ is continuous. }

 However, we are not able to prove this.

 This paper is organized as follows. In Section 2, we deal with the density theorem for packing measure of self-similar sets with SOSC. Firstly, we give the proof of Theorem 1.2 which plays an important role in giving the explicit formula of packing measure in SOSC case. Secondly, we prove the formula $(\ref{05})$ by using the so-called blow-up principle in the SSC case. Section 3 is devoted to the proof of Theorem 1.1 by using the explicit formula $(\ref{05})$.

\section{Density theorems for packing measure of self-similar sets}
We analyze the local behaviour of the packing measure of
self-similar sets in this section. Let $N\geq 2$ be an integer.
$\textbf{f}=\{f_1,f_2,\cdots,f_N\}$ be a IFS on $\mathbb{R}^d$ of
contractive similitudes. In this section, for the sake of
simplicity, we always write $K$ for the self-similar set of
$\textbf{f}$ and we write $s$ for the common value of the Hausdorff
dimension and the packing dimension of $K$. Also,  let $\lambda$
denote the normalised $s$-dimensional Hausdorff measure restricted
to $K$ and write $r_i$ for the contraction ratio of $f_i$   for each
$1\leq i\leq N$ and  write $r_*=\min_{1\leq i\leq N}r_i$. Our main
result in this section, i.e., Theorem 1.2 says that if $K$ satisfies
the SOSC, then
$$\mathcal{P}^s(K\cap B(x,r))\geq (2r)^s$$ for all $B(x,r)\subset \mathcal{O}$ centered in $K$. This result has several applications on densities and can also be applied to compute the exact value of the packing measure $\mathcal{P}^s(K)$ of $K$. Recall that in $\cite{Ols1}$, Olsen  also proved a density theorem for packing measure of self-similar sets which requires that the IFSs satisfy the SSC. In that setting, there exists $r_0>0$ such that the above formula holds for all $x\in K$ and all $r\in(0,r_0]$. However, in the SOSC case, in stead of  finding constant $r_0$, we require $r$ to be small enough such that $B(x,r)\subset \mathcal{O}$. It is easy to check that our result is a natural generalization of the SSC case.

In order to prove Theorem 1.2, we shall need the following lemma.

\textbf{Lemma 2.1.} \textit{Let $B(x,r)\subset \mathcal{O}$ be a
ball centered in $K$ and $k$ a positive integer, then
$$
\mathcal{P}^s(K\cap \bigcup_{\textbf{i}\in
W_k}f_{\textbf{i}}(\overline{B(x,r)}))=\mathcal{P}^s(K\cap\overline{B(x,r)})>0.
$$}

\textit{Proof.} First, we prove that
\begin{equation}\label{2}
K\cap\bigcup_{\textbf{i}\in
W_k}f_\textbf{i}(B(x,r))=\bigcup_{\textbf{i}\in
W_k}f_\textbf{i}(K\cap B(x,r)).
\end{equation}

Fix $y\in K\cap\bigcup_{\textbf{i}\in W_k}f_\textbf{i}(B(x,r))$.
Since $B(x,r)\subset \mathcal{O}$, there exists a $\textbf{u}\in
W_k$ such that $y\in f_\textbf{u}(B(x,r))\subset
f_\textbf{u}(\mathcal{O})$. We also have $y\in
K=\bigcup_{\textbf{i}\in W_k}K_\textbf{i}$ and we therefore find
$\textbf{v}\in W_k$ such that $y\in f_\textbf{v}(K)\subset
f_\textbf{v}(\overline{\mathcal{O}})$. Thus $y\in
f_\textbf{u}(\mathcal{O})\cap f_\textbf{v}(\overline{\mathcal{O}})$,
and therefore $\textbf{u}=\textbf{v}.$ Hence $y\in
f_\textbf{u}(B(x,r))\cap f_\textbf{u}(K)=f_\textbf{u}(K\cap
B(x,r))\subset\bigcup_{\textbf{i}\in W_k}f_\textbf{i}(K\cap
B(x,r)).$ The other direction is obvious. Hence the formula
(\ref{2}) holds.

It follows from (\ref{2}) that
$$\mathcal{P}^s(K\cap\bigcup_{\textbf{i}\in W_k}f_\textbf{i}(B(x,r)))=\mathcal{P}^s(\bigcup_{\textbf{i}\in W_k}f_\textbf{i}(K\cap B(x,r))).$$ However, since the SOSC is satisfied and $B(x,r)\subset \mathcal{O}$ the sets $f_\textbf{i}(K\cap B(x,r))_{\textbf{i}\in W_k}$ are pairwise disjoint. It therefore follows
\begin{eqnarray*}
\mathcal{P}^s(K\cap\bigcup_{\textbf{i}\in
W_k}f_\textbf{i}(B(x,r)))&=&\sum_{\textbf{i}\in
W_k}\mathcal{P}^s(f_\textbf{i}(K\cap
B(x,r)))\\&=&\sum_{\textbf{i}\in W_k}r_\textbf{i}^s
\mathcal{P}^s(K\cap B(x,r))\\&=&\mathcal{P}^s(K\cap B(x,r)).
\end{eqnarray*}

Since the intersection of $K$ with any $n-1$ dimensional $C^1$
manifold  is an $\mathcal{P}^s$-null set. We have
$$\mathcal{P}^s(K\cap\bigcup_{\textbf{i}\in W_k}f_\textbf{i}(B(x,r)))=\mathcal{P}^s(K\cap\bigcup_{\textbf{i}\in W_k}f_\textbf{i}(\overline{B(x,r)})),$$
and
\begin{equation}\label{3}
\mathcal{P}^s(K\cap B(x,r))=\mathcal{P}^s(K\cap \overline{B(x,r)}).
\end{equation}

Using the above three equalities, we get
$$\mathcal{P}^s(K\cap \bigcup_{\textbf{i}\in W_k}f_{\textbf{i}}(\overline{B(x,r)}))=\mathcal{P}^s(K\cap\overline{B(x,r)}).$$
Moreover, since $x\in K$, we deduce that $\mathcal{P}^s(K\cap
\overline{B(x,r)})>0$. This completes the proof of Lemma 1. $\Box$

\textit{Proof of Theorem 1.2.}

In order to reach a contradiction, we assume that $(\ref{1})$ is not
satisfied, i.e., there exists a ball $B(x,r)\subset\mathcal{O}$
centered in $K$, such that
$$\mathcal{P}^s(K\cap B(x,r))<(2r)^s.$$

From $(\ref{3})$, we get
$$\mathcal{P}^s(K\cap \overline{B(x,r)})<(2r)^s.$$
Thus we can find a number $0<\kappa<1$ with
\begin{equation}\label{6}
(1+\kappa)\mathcal{P}^s(K\cap \overline{B(x,r)})<(2r)^s.
\end{equation}

Next, fix $\delta>0$ and choose a positive integer $k$ such that
$$2r_\textbf{i}r\leq \delta$$ for all $\textbf{i}\in W_k$. Let
$\eta=\frac{1}{2}\kappa\mathcal{P}^s(K\cap \overline{B(x,r)}).$ It
follows from Lemma 2.1, $\eta>0.$

For a positive integer $m,$ write
$F_m=K\setminus\bigcup_{\textbf{i}\in
W_k}B(f_\textbf{i}x,r_\textbf{i}r+\frac{1}{m}),$ and observe that
$$F_1\subset F_2\subset F_3\subset\cdots$$
and$$\bigcup_m F_m=K\setminus\bigcup_{\textbf{i}\in
W_k}f_\textbf{i}(\overline{B(x,r)}).$$

If $K\setminus\bigcup_{\textbf{i}\in
W_k}f_\textbf{i}(\overline{B(x,r)})\neq \emptyset,$ then there is a
positive integer ${m_0}$ with $\frac{1}{{m_0}}<\delta$ such that
$F_{m_0}\neq \emptyset$, and
\begin{equation}\label{5}
\mathcal{P}^s(K\setminus\bigcup_{\textbf{i}\in
W_k}B(f_\textbf{i}x,r_\textbf{i}r+\frac{1}{{m_0}}))=\mathcal{P}^s(F_{m_0})\geq
\mathcal{P}^s(K\setminus\bigcup_{\textbf{i}\in
W_k}f_\textbf{i}(\overline{B(x,r)}))-\frac{\eta}{2}.\end{equation}

We can also choose a $\frac{1}{{m_0}}$-packing $\{B(x_i,
\rho_i)\}_i$ of $K\setminus\bigcup_{\textbf{i}\in
W_k}B(f_\textbf{i}x,r_\textbf{i}r+\frac{1}{{m_0}})$ such that
\begin{eqnarray}\label{4}
\nonumber\sum_{i}(2\rho_i)^s&\geq&
P^s_{\frac{1}{{m_0}}}(K\setminus\bigcup_{\textbf{i}\in
W_k}B(f_\textbf{i}x,r_\textbf{i}r+\frac{1}{{m_0}}))-\frac{\eta}{2}\\\nonumber
&\geq& P^s(K\setminus\bigcup_{\textbf{i}\in W_k}B(f_\textbf{i}x,r_\textbf{i}r+\frac{1}{{m_0}}))-\frac{\eta}{2}\\
&\geq& \mathcal{P}^{s}(K\setminus\bigcup_{\textbf{i}\in
W_k}B(f_\textbf{i}x,r_\textbf{i}r+\frac{1}{{m_0}}))-\frac{\eta}{2}.
\end{eqnarray}

Since $x\in K$ and $B(x,r)\subset \mathcal{O}$,
$f_\textbf{i}(B(x,r))\cap f_\textbf{j}(B(x,r))=\emptyset$ for all
$\textbf{i}\neq \textbf{j}$ in $W_k$, and for each $\textbf{i}\in
W_k$, we have $f_\textbf{i}(x)\in K_\textbf{i}\subset K$ and
$2r_\textbf{i}r\leq \delta$. Thus the family
$\{f_\textbf{i}(B(x,r))\}_{\textbf{i}\in W_k}$ is a $\delta$-packing
of $K\cap\bigcup_{\textbf{i}\in
W_k}f_\textbf{i}(\overline{B(x,r)})$.

Since $\{B(x_i,\rho_i)\}_i$ is also  a $\frac{1}{{m_0}}$-packing of
$K\setminus\bigcup_{\textbf{i}\in
W_k}B(f_\textbf{i}x,r_\textbf{i}r+\frac{1}{{m_0}})$, we conclude
that $\{f_\textbf{i}(B(x,r))\}_{\textbf{i}\in
W_k}\bigcup\{B(x_i,\rho_i)\}_i$ is a $\delta$-packing of $K$. Using
this we therefore conclude from $(\ref{6})$, $(\ref{5})$,
$(\ref{4})$, and Lemma 2.1 that
\begin{eqnarray*}
P^s_\delta(K)&\geq& \sum_{\textbf{i}\in W_k}(2r_\textbf{i}r)^s+\sum_i(2\rho_i)^s\\
&\geq& \sum_{\textbf{i}\in W_k}r_\textbf{i}^s(2r)^s+\mathcal{P}^{s}(K\setminus\bigcup_{\textbf{i}\in W_k}B(f_\textbf{i}x,r_\textbf{i}r+\frac{1}{{m_0}}))-\frac{\eta}{2}\\
&\geq& (2r)^s+\mathcal{P}^s(K\setminus\bigcup_{\textbf{i}\in W_k}f_\textbf{i}(\overline{B(x,r)}))-\eta\\
&\geq &(1+\kappa)\mathcal{P}^s(K\cap\overline{B(x,r)})+\mathcal{P}^s(K\setminus\bigcup_{\textbf{i}\in W_k}f_\textbf{i}(\overline{B(x,r)}))-\eta\\
&=& (1+\kappa)\mathcal{P}^s(K\cap\bigcup_{\textbf{i}\in W_k}f_\textbf{i}(\overline{B(x,r)}))+\mathcal{P}^s(K\setminus\bigcup_{\textbf{i}\in W_k}f_\textbf{i}(\overline{B(x,r)}))-\eta\\
&=& \mathcal{P}^s(K)+\kappa\mathcal{P}^s(K\cap\bigcup_{\textbf{i}\in W_k}f_\textbf{i}(\overline{B(x,r)}))-\eta\\
&=& \mathcal{P}^s(K)+2\eta-\eta\\
&=& \mathcal{P}^s(K)+\frac{1}{2}\kappa \mathcal{P}^s(K\cap
\overline{B(x,r)}).
\end{eqnarray*}

Finally, let $\delta\rightarrow 0$, we get
\begin{equation}\label{7}
P^s(K)\geq \mathcal{P}^s(K)+\frac{1}{2}\kappa\mathcal{P}^s(K\cap
\overline{B(x,r)}).
\end{equation}

In $\cite{Fen}$ it is proved that the packing premeasure $P^s$
coincides with the packing measure $\mathcal{P}^s$ for compact
subsets with finite $P^s$-measure. Thus they coincide for $K$, and
it follows from $(\ref{7})$ that
$$\mathcal{P}^s(K)\geq \mathcal{P}^s(K)+\frac{1}{2}\kappa\mathcal{P}^s(K\cap \overline{B(x,r)}).$$
Since $\mathcal{P}^s(K)$ is positive and finite, and
$\frac{1}{2}\kappa\mathcal{P}^s(K\cap \overline{B(x,r)})>0$, we get
the contradiction.

On the other hand, if $K\setminus\bigcup_{\textbf{i}\in
W_k}f_\textbf{i}(\overline{B(x,r)})=\emptyset,$ i.e.,
$K\subset\bigcup_{\textbf{i}\in
W_k}f_\textbf{i}(\overline{B(x,r)})$, then the foregoing string of
inequalities simplifies to
\begin{eqnarray*}
P^s_\delta(K)&\geq& \sum_{\textbf{i}\in W_k}(2r_\textbf{i}r)^s\\
&=& \sum_{\textbf{i}\in W_k}r_\textbf{i}^s(2r)^s\\
&=& (2r)^s\\
&>&(1+\kappa)\mathcal{P}^s(K\cap\overline{B(x,r)})\\
&=& (1+\kappa)\mathcal{P}^s(K\cap\bigcup_{\textbf{i}\in
W_k}f_\textbf{i}(\overline{B(x,r)})).
\end{eqnarray*}
Letting $\delta\rightarrow 0$ and using the fact that
$\mathcal{P}^s(K)=P^s(K)$ this gives
\begin{eqnarray*}
\mathcal{P}^s(K)&=&P^s(K)\\&\geq&  (1+\kappa)\mathcal{P}^s(K\cap\bigcup_{\textbf{i}\in W_k}f_\textbf{i}(\overline{B(x,r)}))\\
&=&(1+\kappa)\mathcal{P}^s(K)>\mathcal{P}^s(K).
\end{eqnarray*}
This provides the desired contradiction.

The proof of Theorem 1.2. is completed. $\Box$

This result has applications on densities. For a given measure $\mu$
on $\mathbb{R}^d$ and $x\in \mathbb{R}^d$, the lower
$\alpha$-density of $\mu$ at $x$ is defined by
$$\Theta_*^\alpha(\mu,x)=\liminf_{r\rightarrow 0}\frac{\mu(B(x,r))}{(2r)^\alpha}.$$
The upper $\alpha$-density $\Theta^{*\alpha}(\mu,x)$ is defined
similarly by taking the upper limit. We have the following result.
If $E\subset \mathbb{R}^d$ and $\alpha>0$ with
$0<\mathcal{P}^\alpha(E)<\infty$, then
\begin{equation}\label{8}
\Theta_*^\alpha(\mathcal{P}^\alpha|_E,x)=1 \mbox{ for }
\mathcal{P}^\alpha-a.e. \quad x\in E.
\end{equation}
See the proof in \cite{Mat1}. We then could get the following
corollary on the basis of $(\ref{8})$ and Theorem 1.2.

\textbf{Corollary 2.2.} \textit{Let $K\subset \mathbb{R}^d$ be the
self-similar set satisfying the SOSC for an open set $\mathcal{O}$,
with packing and Hausdorff dimension $s$. Then
\begin{equation}\label{9}
\mathcal{P}^s(K)=\sup\{\frac{(2r)^s}{\lambda(B(x,r))}| x\in K,
B(x,r)\subset\mathcal{O}\}.
\end{equation}
}

\textit{Proof.} Since $K\cap \mathcal{O}\neq\emptyset,$ we can take
a point $y\in K\cap \mathcal{O}$. Choose $\rho>0$ such that the ball
$B(y,\rho)$ contained in $\mathcal{O}$ and $\mathcal{P}^s(K\cap
B(y,\rho))>0.$ Hence from $(\ref{8})$, there exists a point $z\in
K\cap B(y,\rho)$ with $\Theta_*^s(\mathcal{P}^s|_K,z)=1.$ By the
definition of $\Theta_*^s(\mathcal{P}^s|_K,z)$, there exists a
sequence $\{r_n\}$ with each $r_n\leq\rho-d(z,y)$ and
$r_n\rightarrow 0 $ as $n\rightarrow\infty,$ such that
$\lim_{n\rightarrow\infty}\frac{\mathcal{P}^s(K\cap
B(z,r_n))}{(2r_n)^s}=1.$ Notice that here all balls $B(z,r_n)$ are
contained in $B(y,\rho)\subset\mathcal{O}$. Moreover, by Theorem
1.2, for each ball $B(x,r)\subset\mathcal{O}$ centered in $K$, we
have $\frac{\mathcal{P}^s(K\cap B(x,r))}{(2r)^s}\geq 1$. Hence we
get
$$\inf\{\frac{\mathcal{P}^s(K\cap B(x,r))}{(2r)^s}| x\in K, B(x,r)\subset\mathcal{O}\}=1.$$ Since $\lambda=\frac{\mathcal{P}^s|_K}{\mathcal{P}^s(K)}$, $(\ref{9})$ follows immediately from the above equation. $\Box$

After this work was completed, we learned that Mor\'{a}n
$\cite{Mor1}$ had proved, independently, the same result as
Corollary 2.2. However, his proof is quite different of ours. In
fact, in $\cite{Mor1}$, the so-called self-similar tiling principle
plays a central role in the proof. This principle says that any open
subset $U$ of $K$ can be tiled by a countable set of similar copies
of an arbitrarily given closed set with positive Hausdorff or
packing measure while the tiling is exact in the sense that the part
of $U$ which cannot be covered by the tiles is of null measure. The
continuity theorem is not studied in his paper.

The following lemma will be used in the following corollaries.

\textbf{Lemma 2.3.} (blow-up principle)  \textit{Let $B(x,r)\subset
\mathcal{O}$ centered in $K$. Then for any $1\leq j\leq N$,
$f_j(B(x,r))$ has the same reciprocal density as $B(x,r)$, i.e.,
$\frac{(2r)^s}{\lambda(B(x,r))}=\frac{(2rr_j)^s}{\lambda(f_j(B(x,r)))}$.
In other words, if $B(x,r)\subset f_j(\mathcal{O})$ centered in
$K_j$ for some $1\leq j\leq N,$ then $f_j^{-1}(B(x,r))$ has the same
reciprocal density as $B(x,r).$}

\textit{Proof.} We only need to check that

$$\lambda(f_j(B(x,r)))=r_j^s\lambda(B(x,r)).$$
Actually, $\lambda(f_j(B(x,r)))=\lambda(f_j(B(x,r))\cap
K_j)+\sum_{i\neq j}\lambda(f_j(B(x,r))\cap K_i).$ Notice that if
$i\neq j$, then $f_j(B(x,r))\cap K_i\subset f_j(\mathcal{O})\cap
f_i(K)\subset f_{j}(\mathcal{O})\cap
f_i(\overline{\mathcal{O}})=\emptyset.$ Hence
$\lambda(f_j(B(x,r)))=\lambda(f_j(B(x,r)\cap
K))=r_j^s\lambda(B(x,r)).$ $\Box$

Combining the above lemma and Corollary 2.2, we immediately get the
following corollaries.

\textbf{Corollary 2.4.} \textit{Let $K\subset \mathbb{R}^d$ be the
self-similar set satisfying the SOSC for an open set $\mathcal{O}$,
with packing and Hausdorff dimension $s$. Then
$$
\mathcal{P}^s(K)=\sup\{\frac{(2r)^s}{\lambda(B(x,r))}| x\in K,
B(x,r)\subset\mathcal{O}, B(x,r)\nsubseteq f_j(\mathcal{O}), 1\leq
j\leq N\}.
$$}

\textbf{Corollary 2.5.}  \textit{If $\Delta>0$ and
$d(K_i,K_j)>\Delta$ for all $i\neq j$. Then
$$
\mathcal{P}^s(K)=\sup_{x\in K, \frac{1}{2}r_*\Delta\leq r\leq
\frac{1}{2}\Delta}\frac{(2r)^s}{\lambda(B(x,r))},
$$
where $r_*=\min_{1\leq i\leq N} r_i.$}

\textit{Proof.} Let $\mathcal{O}=\bigcup_{x\in
K}B(x,\frac{\Delta}{2}).$ Then it is obvious that this open set
$\mathcal{O}$ satisfies the SOSC, and therefore the results in
previous can be fully applied in the SSC case. Hence from Corollary
2.2, we get
\begin{equation}\label{40}
\mathcal{P}^s(K)=\sup_{x\in K, 0< r\leq
\frac{1}{2}\Delta}\frac{(2r)^s}{\lambda(B(x,r))}.
\end{equation}

By Lemma 2.3, we can limit $B(x,r)$ not contained in each
$f_j(\mathcal{O})$. Fix $x\in K$, $0<r\leq \frac{1}{2}\Delta$. Then
$B(x,r)\subset\mathcal{O}$ is a ball centered in $K$. Hence there
exists $j$, such that $x\in K_j$. Obviously $B(x,r)\cap K_j\neq
\emptyset,$ which yields that $B(x,r)\cap
f_j(\overline{\mathcal{O}})\neq\emptyset$. Hence $B(x,r)\cap
f_j({\mathcal{O}})\neq\emptyset$. But $B(x,r)$ can not contained in
$f_j(\mathcal{O})$, i.e., $B(x,r)\nsubseteq\bigcup_{z\in
{K_j}}{B(z,\frac{r_j}{2}\Delta)}$. So $r\geq \frac{r_j}{2}\Delta\geq
\frac{r_*}{2}\Delta.$ Hence
$$\mathcal{P}^s(K)=\sup_{x\in K, \frac{1}{2}r_*\Delta\leq r\leq \frac{1}{2}\Delta}\frac{(2r)^s}{\lambda(B(x,r))}. \Box$$

\section{Proof of Theorem 1.1.}
In this section, we prove Theorem 1.1. In order to prove this
theorem, we need some lemmas. Below $D_H(\cdot,\cdot)$ denotes the
Hausdorff metric on the family of all compact subsets of $X$.

\textbf{Lemma 3.1.} \textit{Let $\Delta>0$, $\textbf{f}\in
M_\Delta$. Then there exists $\beta>0$ such that if $\textbf{g}\in
M_\Delta$ with $D(\textbf{f},\textbf{g})\leq \beta$, then
$r_*(\textbf{g})\geq \frac{1}{2}r_*(\textbf{f})$.}

\textit{Proof.} It is easy since the map $\textbf{g}\rightarrow
r_*(\textbf{g})$ from $M_\Delta$ into $\mathbb{R}$ is continuous.
$\Box$

\textbf{Lemma 3.2.} \textit{Let $\Delta>0$, $\textbf{f}\in
M_\Delta$, $0<\rho<\frac{1}{4}r_*(\textbf{f})\Delta,$ $\gamma>0$ and
let $\beta$ be the same as that in Lemma 3.1. Then there exists
$\delta_{1}>0$ with $\delta_{1}\leq\beta$ such that if
$\textbf{g}\in M_\Delta$ and $\textbf{h}\in M_\Delta$ with
$D(\textbf{f},\textbf{g})\leq \delta_{1}$ and
$D(\textbf{f},\textbf{h})\leq\delta_{1}$, then for each ball
$B(x,r)$ centered in $K(\textbf{g})$ with radius
$r\in[\frac{1}{2}r_*(\textbf{g})\Delta,\frac{1}{2}\Delta]$, there
exists a ball $B(y,r-\rho)$ centered in $K(\textbf{h})$ such that
$$
\lambda(\textbf{h})(B(y,r-\rho))-\gamma\leq
\lambda(\textbf{g})(B(x,r)).
$$}

\textit{Proof.} Choose a positive integer $k$ such that
$$\mbox{diam} (f_\textbf{i} (K(\textbf{f})))\leq \frac{\rho}{8}.$$
for all $\textbf{i}\in W_k$. By the continuity of  the map
$\textbf{g}\rightarrow K(\textbf{g})$ from $M_\Delta$ into
$\mathcal{K}$ and the map $\textbf{g}\rightarrow r_i(\textbf{g})$
from $M_\Delta$ into $\mathbb{R}$ for each $1\leq i\leq N$, we can
choose $\delta_{1}>0$ with $\delta_{1}\leq \beta$ such that if
$\textbf{g}\in M_\Delta$ with
$D(\textbf{f},\textbf{g})\leq\delta_{1}$, then
\begin{equation}\label{13}
\mbox{diam}( g_{\textbf{i}}(K(\textbf{g})))\leq\frac{\rho}{4},
\end{equation}
\begin{equation}\label{14}
D_H(f_\textbf{i}(K(\textbf{f})),
g_{\textbf{i}}(K(\textbf{g})))\leq\frac{\rho}{8},
\end{equation}
\begin{equation}\label{15}
|(r_\textbf{i}(\textbf{f}))^{s(\textbf{f})}-(r_\textbf{i}(\textbf{g}))^{s(\textbf{g})}|\leq\frac{\gamma}{2N^k}
\end{equation}
for all $\textbf{i}\in W_k$.

Now fix $\textbf{g},\textbf{h}\in M_\Delta$ with
$D(\textbf{f},\textbf{g})\leq\delta_{1}$,
$D(\textbf{f},\textbf{h})\leq\delta_{1}$, also fix a ball $B(x,r)$
centered in $K(\textbf{g})$ with radius
$r\in[\frac{1}{2}r_*(\textbf{g})\Delta,\frac{1}{2}\Delta]$. By Lemma
3.1, $r_*(\textbf{g})\geq \frac{1}{2}r_*(\textbf{f})$, so
$r-\rho>0.$

Since $D_H(K(\textbf{f}),K(\textbf{g}))\leq \frac{\rho}{8}$,
$D_H(K(\textbf{f}),K(\textbf{h}))\leq \frac{\rho}{8}$, we get
$D_H(K(\textbf{g}),K(\textbf{h}))< \frac{\rho}{2}$. Hence there
exists a point $y\in B(x,\frac{\rho}{2})$ with $y\in K(\textbf{h}),$
which yields that
\begin{equation}\label{16}
B(y,r-\rho)\subset B(x,r-\frac{\rho}{2}).
\end{equation}

If we denote $V_k=\{\textbf{i}| \textbf{i}\in W_k,
h_\textbf{i}(K(\textbf{h}))\cap B(x,r-\frac{\rho}{2})\neq
\emptyset\}$, then from $(\ref{15})$ and $(\ref{16})$ we get that
\begin{eqnarray*}
\lambda(\textbf{h})(B(y,r-\rho))-\gamma&\leq&\lambda(\textbf{h})(B(x,r-\frac{\rho}{2}))-\gamma\\
&\leq&\lambda(\textbf{h})(\bigcup_{\textbf{i}\in V_k}h_\textbf{i}(K(\textbf{h})))-\gamma\\
&=&\sum_{\textbf{i}\in V_k}\lambda(\textbf{h})(h_\textbf{i}(K(\textbf{h})))-\gamma\\
&=&\sum_{\textbf{i}\in V_k}(r_\textbf{i}(\textbf{h}))^{s(\textbf{h})}-\gamma\\
&\leq&\sum_{\textbf{i}\in V_k}((r_\textbf{i}(\textbf{f}))^{s(\textbf{f})}+\frac{\gamma}{2N^k})-\gamma\\
&\leq&\sum_{\textbf{i}\in V_k}((r_\textbf{i}(\textbf{g}))^{s(\textbf{g})}+\frac{\gamma}{N^k})-\gamma\\
&\leq&\sum_{\textbf{i}\in V_k}(r_\textbf{i}(\textbf{g}))^{s(\textbf{g})}\\
&=&\sum_{\textbf{i}\in
V_k}\lambda(\textbf{g})(g_\textbf{i}(K(\textbf{g}))).
\end{eqnarray*}

Next, if we denote by $U_k=\{\textbf{i}| \textbf{i}\in W_k,
g_\textbf{i}(K(\textbf{g}))\cap B(x,r-\frac{\rho}{4})\neq
\emptyset\}$, then we must have $V_k\subset U_k$. In fact, if
$\textbf{i}\in V_k$, then $h_\textbf{i}(K(\textbf{h}))\cap
B(x,r-\frac{\rho}{2})\neq \emptyset$. Combining this with
$(\ref{14})$ we get $f_\textbf{i}(K(\textbf{f}))\cap
B(x,r-\frac{3\rho}{8})\neq \emptyset$, and using $(\ref{14})$ once
more we could get that $g_\textbf{i}(K(\textbf{g}))\cap
B(x,r-\frac{\rho}{4})\neq \emptyset$ which proves $V_k\subset U_k.$
Hence,
$$
\lambda(\textbf{h})(B(y,r-\rho))-\gamma\leq\sum_{\textbf{i}\in
U_k}\lambda(\textbf{g})(g_\textbf{i}(K(\textbf{g})))=\lambda(\textbf{g})(\bigcup_{\textbf{i}\in
U_k}g_\textbf{i}(K(\textbf{g}))).
$$

Finally, notice that if $\textbf{i}\in W_k$ and
$g_\textbf{i}(K(\textbf{g}))\cap B(x,r-\frac{\rho}{4})\neq
\emptyset$, then it follows from $(\ref{13})$ that
$g_\textbf{i}(K(\textbf{g}))\subset B(x,r)$, whence
$\bigcup_{\textbf{i}\in U_k}g_\textbf{i}(K(\textbf{g}))\subset
B(x,r)$.  Hence we have
$$\lambda(\textbf{h})(B(y,r-\rho))-\gamma\leq \lambda(\textbf{g})(B(x,r)). \Box$$

\textbf{Lemma 3.3.} \textit{Let $\Delta>0$, $\textbf{f}\in
M_\Delta$, $\kappa>0$, and let $\beta$ be the same as that in Lemma
3.1. Then there exists
$0<\rho<\min\{\beta,\frac{1}{4}r_*(\textbf{f})\Delta\}$ such that if
$\textbf{g},\textbf{h}\in M_\Delta$ with
$D(\textbf{f},\textbf{g})\leq\rho$ and $D(\textbf{f},\textbf{h})\leq
\rho$, and
$r\in[\frac{1}{2}r_*(\textbf{g})\Delta,\frac{1}{2}\Delta]$, then
$$
(2(r-\rho))^{s(\textbf{h})}\geq (2r)^{s(\textbf{g})}-\kappa.
$$}

\textit{Proof.} Notice that the map $(r,\textbf{g})\rightarrow
(2r)^{s(\textbf{g})}$ from $(0,\frac{1}{2}\Delta)\times M_\Delta$
into $\mathbb{R}$ is continuous. Moreover, for each fixed
$\textbf{g}$, the map is uniformly continuous on
$r\in(0,\frac{1}{2}\Delta)$. Hence, for fixed $\textbf{f}\in
M_\Delta$, $\kappa>0$, there exists
$0<\rho<\min\{\beta,\frac{1}{4}r_*(\textbf{f})\Delta\}$ with the
following property. If $\textbf{g}\in M_\Delta$ with
$D(\textbf{f},\textbf{g})\leq\rho$, and
$r',r''\in(0,\frac{1}{2}\Delta)$ with $|r'-r''|\leq\rho$, then
\begin{equation}\label{18}
|(2r')^{s(\textbf{g})}-(2r'')^{s(\textbf{f})}|\leq\frac{\kappa}{2}.
\end{equation}

Hence if we take $\textbf{g}\in M_\Delta$ with
$D(\textbf{f},\textbf{g})\leq\rho$, and $r'=r, r''=r-\rho$ with
$r\in[\frac{1}{2}r_*(\textbf{g})\Delta,\frac{1}{2}\Delta]$, then by
$(\ref{18})$, we get
\begin{equation*}
|(2r)^{s(\textbf{g})}-(2(r-\rho))^{s(\textbf{f})}|\leq\frac{\kappa}{2}.
\end{equation*}

And if we take $\textbf{h}\in M_\Delta$ with
$D(\textbf{f},\textbf{h})\leq\rho$, and $r'=r''=r-\rho$ with
$r\in[\frac{1}{4}r_*(\textbf{f})\Delta,\frac{1}{2}\Delta]$, then
also by $(\ref{18})$, we get
\begin{equation*}
|(2(r-\rho))^{s(\textbf{h})}-(2(r-\rho))^{s(\textbf{f})}|\leq\frac{\kappa}{2}.
\end{equation*}

The above two inequalities and Lemma 3.1 gives the desired result.
$\Box$

\textbf{Lemma 3.4.} \textit{Let $c,$ $C$, $\kappa>0$,
$0<\varepsilon<1$ with $\kappa<c\varepsilon$ and $c<C$. Then there
exists $\gamma>0$
 such that
 $$\frac{x-\kappa}{y+\gamma}\geq \frac{x}{y}-\varepsilon$$
 for all $x,y\in [c,C].$}

\textit{ Proof.} Take $\gamma=\frac{c( c\varepsilon-\kappa)}{C}$.
Without losing generality, we may assume that
$\frac{x}{y}-\varepsilon>0$, i.e., $x-y\varepsilon>0$, then
$\gamma\leq \frac{y(\varepsilon y-\kappa)}{x-y\varepsilon}$. Hence
$(x-y\varepsilon)\gamma\leq y(\varepsilon y-\kappa)$, which yields
$\kappa y+x\gamma\leq \varepsilon y(y+\gamma)$. Thus, dividing the
above inequality by $y(y+\gamma)>0$ gives $\frac{\kappa
y+x\gamma}{y(y+\gamma)}\leq \varepsilon$, i.e.,
$\frac{x}{y}-\frac{x-\kappa}{y+\gamma}\leq \varepsilon.$ $\Box$

\newpage

\textit{Proof of Theorem 1.1.}

Since $M_{SSC}=\bigcup_{\Delta>0}M_\Delta$ and that $M_\Delta$ are
open subsets of $M_{SSC}$ for all $\Delta>0$, we only need to prove
for each $\Delta>0$, the map $$\textbf{f}\rightarrow
\mathcal{P}^{s(\textbf{f})}(K(\textbf{f}))$$ from $M_\Delta$ into
$\mathbb{R}$ is continuous.

Fix $\Delta>0$, $\textbf{f}\in M_\Delta$ and let $0<\varepsilon<1$.
We now find $\delta>0$ such that if $\textbf{g}\in M_\Delta$ with
$D(\textbf{f},\textbf{g})\leq \delta$, then
$$|\mathcal{P}^{s(\textbf{f})}(K(\textbf{f}))-\mathcal{P}^{s(\textbf{g})}(K(\textbf{g}))|\leq
\varepsilon.$$

By the continuity of the map $\textbf{g}\rightarrow s(\textbf{g})$
from $M_\Delta$ into $\mathbb{R}$, there exists $\delta_2>0$ with
$\delta_2\leq\beta$ such that if $\textbf{g}\in M_\Delta$ with
$D(\textbf{f},\textbf{g})\leq\delta_2$, then
$$\frac{1}{2}s(\textbf{f})\leq s(\textbf{g})\leq\frac{3}{2}s(\textbf{f}),$$ where $\beta$ is the same as that in Lemma 3.1.
Hence there exists $C_1, C_2>0$ such that if $\textbf{g}\in
M_\Delta$ with $D(\textbf{f},\textbf{g})\leq\delta_2$, and
$r\in[\frac{1}{2}r_*(\textbf{g})\Delta,\frac{1}{2}\Delta]$, then
$$C_1\leq (2r)^{s(\textbf{g})}\leq C_2.$$

In Lemma 3.1 of \cite{Ols}, the map
$\textbf{g}\rightarrow\lambda(\textbf{g})$ from $M_\Delta$ into
$\mathcal{M}$ is continuous, where $\mathcal{M}$ denotes the space
consist of all Borel regular probability measures equipped with the
weak topology. Hence there exists $\delta_3>0$, $C_3>0$ with
$\delta_3\leq\beta $ such that if $\textbf{g}\in M_\Delta$ with
$D(\textbf{f},\textbf{g})\leq \delta_3$, then
$$\lambda(\textbf{g})(B(x,r))\geq C_3$$ for all ball $B(x,r)$ with radius $r\in[\frac{1}{2}r_*(\textbf{g})\Delta,\frac{1}{2}\Delta]$ centered in $K(\textbf{g}).$

Put $c=\min\{C_1,C_3\}$, $C=C_2+1+c$,
$\kappa=\frac{c\varepsilon}{4}$, $\gamma>0$ as in Lemma 3.4. Take
$\delta=\min\{\rho,\delta_1,\delta_2,\delta_3\}$. We claim that if
$\textbf{g}\in M_\Delta$ with $D(\textbf{f},\textbf{g})\leq\delta$,
then
$$|\mathcal{P}^{s(\textbf{f})}(K(\textbf{f}))-\mathcal{P}^{s(\textbf{g})}(K(\textbf{g}))|\leq \varepsilon.$$
To prove this we show that if $\textbf{g},\textbf{h}\in M_\Delta$
with $D(\textbf{f},\textbf{g})\leq\delta$ and
$D(\textbf{f},\textbf{h})\leq\delta$, then
\begin{equation}\label{30}\mathcal{P}^{s(\textbf{h})}(K(\textbf{h}))\geq\mathcal{P}^{s(\textbf{g})}(K(\textbf{g}))-\varepsilon.\end{equation}

We therefore fix $\textbf{g},\textbf{h}\in M_\Delta$ satisfying
$D(\textbf{f},\textbf{g})\leq\delta$ and
$D(\textbf{f},\textbf{h})\leq\delta$. It follows from the Corollary
2.5 that there exists $B(x,r)$ centered in $K(\textbf{g})$ with
radius $r\in[\frac{1}{2}r_*(\textbf{g})\Delta,\frac{1}{2}\Delta]$
such that
\begin{equation}\label{20}
\frac{(2r)^{s(\textbf{g})}}{\lambda(\textbf{g})(B(x,r))}\geq
\mathcal{P}^{s(\textbf{g})}(K(\textbf{g}))-\frac{\varepsilon}{2}.
\end{equation}

Since $c\leq (2r)^{s(\textbf{g})}\leq C$, $c\leq
\lambda(\textbf{g})(B(x,r))\leq C$ and
$\kappa=\frac{c\varepsilon}{4}<\frac{c\varepsilon}{2}$, then from
Lemma 3.4,
\begin{equation}\label{21}
\frac{(2r)^{s(\textbf{g})}-\kappa}{\lambda(\textbf{g})(B(x,r))+\gamma}\geq
\frac{(2r)^{s(\textbf{g})}}{\lambda(\textbf{g})(B(x,r))}-\frac{\varepsilon}{2}.
\end{equation}

And from Lemma 3.3, we have
\begin{equation}\label{22}
(2(r-\rho))^{s(\textbf{h})}\geq (2r)^{s(\textbf{g})}-\kappa.
\end{equation}

Then from Lemma 3.2, there exists a ball $B(y,r-\rho)$ centered in
$K(\textbf{h})$ such that
\begin{equation}\label{24}
\lambda(\textbf{h})(B(y,r-\rho))-\gamma\leq
\lambda(\textbf{g})(B(x,r)).
\end{equation}

Combining $(\ref{40})$ and $(\ref{20})$ to $(\ref{24})$, we get
\begin{eqnarray*}
\mathcal{P}^{s(\textbf{h})}(K(\textbf{h}))&\geq&\frac{2(r-\rho)^{s(\textbf{h})}}{\lambda(\textbf{h})(B(y,r-\rho))}\\
&\geq&\frac{(2r)^{s(\textbf{g})}-\kappa}{\lambda(\textbf{g})(B(x,r))+\gamma}\\
&\geq&\frac{(2r)^{s(\textbf{g})}}{\lambda(\textbf{g})(B(x,r))}-\frac{\varepsilon}{2}\\
&\geq& \mathcal{P}^{s(\textbf{g})}(K(\textbf{g}))-\varepsilon.
\end{eqnarray*}
This proves $(\ref{30})$ and hence the proof of Theorem 1.1 is
completed. $\Box$

\textbf{Acknowledgements} I would like to thank the referee for his
valuable comments and suggestions, especially for the refinement of
the proof of Theorem 1.2, that led to the improvement of the
manuscript. This project has been supported by the National Natural
Science Foundation of China 10901081.

\end{document}